\documentclass[11pt,a4paper]{article}

\DeclareMathAlphabet{\mathcal}{OMS}{cmsy}{m}{n}
\usepackage[top=25mm,bottom=35mm,left=20mm,right=20mm]{geometry} 

\usepackage{url}
\usepackage{graphicx}
\usepackage{mathptmx}     
\usepackage{mathtools}
\mathtoolsset{showonlyrefs=true}
\usepackage{color}
\usepackage{doi}
\usepackage{amsmath}
\usepackage{amsfonts}
\usepackage{amssymb}
\usepackage{amsthm}
\usepackage{bm}

\usepackage{algorithm}
\usepackage{algorithmic}
\usepackage{bm}
\newcommand{\bx}{\bm{x}}
\newcommand{\by}{\bm{y}}
\newcommand{\bb}{\bm{b}}
\newcommand{\br}{\bm{r}}
\newcommand{\bc}{\bm{c}}
\newcommand{\bp}{\bm{p}}

\newcommand{\bbR}{\mathbb{R}}
\newcommand{\rmd}{\mathrm{d}}
\newcommand{\dg}{\nabla _\rmd}

\DeclareMathOperator*{\argmin}{arg\,min}
\DeclareMathOperator*{\diag}{diag}

\newcommand{\Go}{G_\omega}
\newcommand{\co}{\bc_\omega}

\theoremstyle{definition}
\newtheorem{proposition}{Proposition}

\newtheorem{theorem}{Theorem}
\newtheorem{remark}{Remark}

\usepackage{tikz}
\usepackage{pgfplotstable}
\usetikzlibrary{arrows,positioning,plotmarks,external,patterns,angles,
decorations.pathmorphing,backgrounds,fit,shapes,graphs,calc}

\begin{document}
\title{Adaptive SOR methods based on the Wolfe conditions}
\author{
Yuto Miyatake\thanks{Cybermedia Center, Osaka University, 
			1-32 Machikaneyama, Toyonaka, Osaka 560-0043, Japan,
\href{mailto:miyatake@cas.cmc.osaka-u.ac.jp}{miyatake@cas.cmc.osaka-u.ac.jp}} ,\
Tomohiro Sogabe\thanks{Department of Applied Physics, 
              Graduate School of Engineering, 
              Nagoya University,  Furo-cho, Chikusa-ku, 
              464-8603 Nagoya, Japan, 
\href{mailto:sogabe@na.nuap.nagoya-u.ac.jp}{sogabe@na.nuap.nagoya-u.ac.jp}}
\ and 
Shao-Liang Zhang\thanks{Department of Applied Physics, 
              Graduate School of Engineering, 
              Nagoya University,  Furo-cho, Chikusa-ku, 
              464-8603 Nagoya, Japan, 
\href{mailto:zhang@na.nuap.nagoya-u.ac.jp}{zhang@na.nuap.nagoya-u.ac.jp}}
}
%\date{August, 2015}

\maketitle

\begin{abstract}
Because the expense of estimating the optimal value of the relaxation parameter in the successive over-relaxation (SOR) method
is usually prohibitive, the parameter is often adaptively controlled.
In this paper, new adaptive SOR methods are presented that
are applicable to a variety of symmetric positive definite linear systems and do not require additional matrix-vector products when updating the parameter.
To this end, we regard the SOR method as an algorithm for minimising a certain objective function, which yields an interpretation of the relaxation parameter as the step size
following a certain change of variables.
This interpretation enables us to adaptively control the step size based on some line search techniques,
such as the  Wolfe conditions.
Numerical examples demonstrate the favourable behaviour of the proposed methods.
\end{abstract}

\section{Introduction}
\label{sec1}
The successive over-relaxation (SOR) method is 
one of the most well-known stationary iterative methods for solving linear systems
\begin{align}
A\bx = \bb ,
\end{align}
where $A \in \bbR^{n\times n}$ and $ \bb \in \bbR^{n}$ are given nonsingular matrix and 
vector, respectively.
The SOR method is formulated as
\begin{align}
\bx^{(k+1)} = \Go \bx^{(k)} + \co, \quad k = 0,1,2,\dots,
\end{align}
where
\begin{align}
\Go &= (D+\omega L )^{-1} [ (1-\omega) D - \omega U], \\
\co &= \omega (D+\omega L)^{-1} \bb.
\end{align}
Here, the matrix $A$ is expressed as the matrix sum
\begin{align}
A = D+L+U,
\end{align}
where $D = \diag (a_{11}, a_{22},\dots, a_{nn})$, and $L$ and $U$ are
strictly lower and upper triangular $n\times n $ matrices.

The convergence performance of the SOR method depends on the relaxation parameter $\omega$.
A necessary condition for the SOR method to converge is that $\omega \in (0,2)$,
and this condition is also sufficient for a symmetric positive definite matrix $A$
(for further details of SOR theory see, e.g.,~\cite{ha00,va00}).
The relaxation parameter $\omega$
is ideally chosen such that the spectral radius of the iteration matrix $\Go$ is 
as small as possible.
However, the optimal value is rarely available except in some
special cases, such as consistently ordered $p$-cyclic matrices~\cite{va00}.
Thus, in practice a specific value obtained either empirically or heuristically is 
often used. 
A typical example is $\omega = 2-O(\tilde{h})$, where $\tilde{h}$ 
is the mesh spacing of the 
discretisation of the underlying physical domain.
An alternative approach is to adaptively control the relaxation parameter.
This approach is referred to as an adaptive SOR method.

In this paper, we are concerned with adaptive SOR methods.
To explain our motivation, we first review two standard approaches and their features.
One is to approximate the spectral radius of $\Go$ using the information obtained during iterations~\cite{hy81}.
This approach is based on the observation that
for a consistently ordered coefficient matrix with some additional assumptions,
the spectral radius of $\Go$ is related to that of the Jacobi iteration matrix.
While a highly accurate parameter is often estimated by this method, and
the proportion of the cost of the estimation is small compared with a whole iteration,
the applicability depends on the matrix properties.
The second approach is to approximately minimise the residual instead 
of the spectral radius of $\Go$, in order
to handle more general linear systems~\cite{ba03} (see~\cite{me04,rr16} for
some related approaches).
The key observation in~\cite{ba03} is that the residual 
can be expanded as a polynomial
in terms of $\omega$, and the authors' approach is to
first approximate or truncate the polynomial by a lower order one, 
then 
determine the real root of it.
For symmetric positive definite linear systems,
the authors also considered a certain quadratic form instead of the residual.
While the efficiency was reported, 
the estimated parameter does not always satisfy the convergence condition, and
several additional matrix-vector products are required to update the parameter.

By taking this background into account, the aim of this paper is
to develop a new type of adaptive SOR method that
is applicable to a variety of symmetric positive definite linear systems,
without additional assumptions or a requirement for additional
matrix-vector products when updating the relaxation parameter.
To achieve this goal, we develop new adaptive SOR methods 
based on our previous work~\cite{ms18},
which shows that for any symmetric positive definite linear system
the SOR method can be regarded as an algorithm for solving 
a certain minimisation problem.
The consequence of this discussion is that 
the relaxation parameter can be interpreted as the step size
following the change of variables $h = 2\omega / (2-\omega)$.
This enables us to apply some line search techniques, developed 
in the context of unconditioned optimisation problems.
In~\cite{ms18}, we already tested an approach 
based on the steepest descent method,
but this approach requires an additional matrix-vector product 
to update the relaxation parameter,
and the convergence could be slow for certain linear systems, 
as will be illustrated later.
%In this work, instead of the steepest descent method
%we shall apply the Armijo condition and Wolfe conditions with
%slight modifications, to propose two new adaptive SOR methods
%(one is based on the Armijo condition only, and the other on the Wolfe conditions).
In this work, 
we shall propose two new adaptive SOR methods based on other line search techniques.
One is based on the Armijo condition, which makes the computational cost for updating the parameter inexpensive.
But the estimated parameter sometimes tends to $0$, which delays the convergence. To avoid such a situation, we consider another adaptive SOR method based on the Wolfe conditions, i.e. the curvature condition
is also used in addition to the Armijo condition.
These methods restrict $A$ to be symmetric positive definite, 
but do not require additional
assumptions on the coefficient matrix.
Furthermore, they do not require additional matrix-vector products when updating the parameter, 
and it is guaranteed that the estimated value always satisfies the convergence condition $\omega \in (0,2)$.

The remainder of this paper is organised as follows.
In Section~\ref{sec2}, the connection between the SOR method and an optimisation problem
is explained, based on~\cite{ms18}.
New adaptive SOR methods are presented in Section~\ref{sec3}, and 
numerically 
tested in Section~\ref{sec4}.
Finally, concluding remarks are provided in Section~\ref{sec5}.

\section{Connection between the SOR method and an optimisation problem}
\label{sec2}

This section briefly reviews our previous work~\cite{ms18}, which
reveals the connection between the SOR method and an optimisation problem.

\subsection{Unconditioned optimisation problem for a strictly convex objective function}

Let a differentiable function $f:\bbR^n \to  \bbR$ be 
strictly convex\footnote{A function $f\in\bbR^n\to\bbR$
is said to be strictly convex if and only if 
for all $\bx,\by\in\bbR^n$ ($\bx\neq\by$)
and $\lambda\in(0,1)$, it holds that $f(\lambda \bx + (1-\lambda)\by ) < \lambda f(\bx) 
+ (1-\lambda) f(\by)$.} and coercive\footnote{A function $f\in\bbR^n\to\bbR$
is said to be coercive if and only if $f(\bx) \to \infty$ for $\| \bx\| \to\infty$.}. 
Let us consider the unconditioned optimisation problem
\begin{align} \label{eq:op}
\min_{\bx\in\bbR^n} f(\bx).
\end{align}
Then, there exists a unique global minimiser $\bx_\ast$ for the function $f$:
\begin{align}
\bx_\ast = \argmin_{\bx\in\bbR^n} f(\bx).
\end{align}
Note that
this minimiser coincides with the unique solution to $\nabla f(\bx ) = \bm{0} $,
and can also be realised as the equilibrium of the 
following gradient system:
\begin{align}\label{eq:gs}
\frac{\rmd}{\rmd t}\bx (t) = -P \nabla f(\bx(t)), 
\end{align}
where $P\in\bbR^{n\times n}$ is a symmetric positive definite matrix.
More precisely, the solution to \eqref{eq:gs} converges to 
the unique minimiser of the function $f$
for any initial vector $\bx(0)=\bx_0$, i.e.,
\begin{align*}
\lim _{t\to \infty} \bx (t) = \argmin_{\bx\in\bbR^n} f(\bx).
\end{align*}

In the following, we discuss algorithms for the optimisation problem \eqref{eq:op} 
with an emphasis on the relation to the gradient system \eqref{eq:gs}.
Note that most gradient descent methods 
can be interpreted as the explicit Euler method 
for the gradient system \eqref{eq:gs}:
\begin{align} \label{eq:eem}
\bx ^{(k+1)} = \bx^{(k)} - h P \nabla f(\bx^{(k)}),
\end{align}
where $h$ denotes the step size.
In particular, \eqref{eq:eem} coincides with the steepest descent method 
when $P$ is the identity matrix.

Instead of the explicit Euler method, 
we consider an alternative discretisation method for the gradient system \eqref{eq:gs}.
The key is to approximate the gradient using 
the so-called discrete gradient~\cite{go96,hl14,ia88,mf14,mq99,mqr98,qm08,qt96}.
We consider a discrete approximation to $\nabla f$,
which will be denoted by $\dg f$.
The function $\dg f : \bbR^n \times \bbR^n \to \bbR$ is defined by
\begin{align}\def\arraystretch{1.3}
\dg f(\bx,\by) = 
\begin{bmatrix}
\cfrac{f(x_1,y_2,\dots,y_n)-f(y_1,\dots, y_n)}{x_1-y_1} \\
\cfrac{f(x_1,x_2,y_3\dots,y_n)-f(x_1,y_2,\dots, y_n)}{x_2-y_2} \\
\vdots\\
\cfrac{f(x_1,\dots,x_n)-f(x_1,\dots,x_{n-1}, y_n)}{x_n-y_n}
\end{bmatrix}  \label{eq:ia}
\end{align}
for all $\bx,\by \in \bbR^n$.
This function is called a discrete gradient, 
because it is an approximation to the gradient
in the sense that
\begin{align} \label{eq:co}
\lim_{\by \to \bx} \dg f(\bx,\by ) = \nabla f(\bx).
\end{align}
Furthermore, it follows that
\begin{align} \label{dcr}
f(\bx) - f(\by) = \dg f(\bx,\by)^\top (\bx-\by)
\end{align}
for all $\bx,\by \in \bbR^n$.
By utilising the discrete gradient \eqref{eq:ia},
we discretise the gradient system \eqref{eq:gs} as follows:
\begin{align}
\bx ^{(k+1)} = \bx^{(k)} - hP \dg f(\bx^{(k+1)},\bx^{(k)}).
%\frac{\bx^{(k+1)} - \bx^{(k)} }{h} = - P \dg f(\bx^{(k+1)},\bx^{(k)}),
\label{dgscheme}
\end{align}
Here, we have just replaced the gradient in \eqref{eq:eem} with the discrete gradient.
However, there is a significant difference between \eqref{eq:eem} and \eqref{dgscheme} from
the viewpoint of optimisation problems.
For the solution to~\eqref{dgscheme},
the following discrete dissipation property holds:
\begin{align}
\frac{1}{h} \left( f(\bx^{(k+1)}) - f(\bx^{(k)}) \right)
&=   \dg f(\bx^{(k+1)}, \bx^{(k)} ) ^\top \frac{\bx^{(k+1)} - \bx ^{(k)}}{h}  \\
&= - \dg f(\bx^{(k+1)}, \bx^{(k)} ) ^\top P \dg f(\bx^{(k+1)}, \bx^{(k)} ) \leq 0 \label{d:dp}
\end{align}
as long as $h>0$.
Thanks to this property, the sequence obtained by~\eqref{dgscheme} converges to the minimiser.

\begin{proposition}[see, e.g.,~\cite{gr17}] \label{prop:d}
Let a differentiable function $f:\bbR^n\to\bbR$ be strictly convex and coercive.
Then, the sequence $\{ \bx ^{(k)}\}_{k=0}^\infty$ obtained by the iteration~\eqref{dgscheme}
converges to the unique minimiser of the function $f$ for any initial vector $\bx_0$,
i.e.,
\begin{align}
\lim _{k\to\infty} \bx^{(k)} = \argmin _{\bx\in\bbR^n} f(\bx).
\end{align}
\end{proposition}
A rigorous proof of this is given in~\cite{gr17}.

This proposition indicates that the convergence is guaranteed even for a relatively large step size $h$, which is not the case for standard gradient descent methods \eqref{eq:eem}.

Let us mention the case that the step size $h$ is controlled adaptively.
We denote the step size at $k$th iteration by $h^{(k)}$.
In this case, we need to avoid the situation that the step size tends to $0$ or $+\infty$.
The convergence
is guaranteed if $0 < \epsilon < h^{(k)} < M < \infty$ for two positive constants $\epsilon$ and $M$.

\begin{remark}
In the context of numerical ordinary differential equations, 
a function satisfying both \eqref{eq:co} and \eqref{dcr} is called a discrete gradient,
and is used to simulate systems of ordinary differential equations such as gradient systems and Hamiltonian systems.
Such a function is not unique in general, and \eqref{eq:ia} is just an example, which was
introduced by Itoh--Abe~\cite{ia88}.
See, e.g.,~\cite{go96,hl14,ia88,mf14,mq99,mqr98,qm08,qt96}
for further details regarding the discrete gradient method.
\end{remark}

\subsection{SOR method as an optimisation solver}

We now show that the SOR method 
can be regarded as an algorithm to solve a certain optimisation problem.

Let 
\begin{align} \label{qf}
f(\bx) = \frac12 \bx^\top A \bx  - \bx^\top \bb.
\end{align}
If $A$ is symmetric positive definite, then the function $f$ is strictly convex and coercive,
and thus it yields the unique minimiser 
$\bx_\ast = A^{-1}\bb$ (note that $\nabla f(\bx) = A\bx - \bb$).
In particular,
when $P=D^{-1}$
the scheme \eqref{dgscheme} is written as
\begin{align}\label{dgscheme1}
\bx^{(k+1)} = \bx^{(k)} -hD^{-1} \dg f(\bx^{(k+1)}, \bx^{(k)}), 
\end{align}
where
the $i$th component of the discrete gradient is calculated to be
\begin{align}
(\dg f(\bx,\by))_i &= \frac{f(x_1,\dots,x_i,y_{i+1},\dots,y_n)-f(x_1,\dots,x_{i-1},y_i,\dots, y_n)}{x_i-y_i} \\
&= \sum_{j<i} a_{ij}x_j + a_{ii}\frac{x_i+y_i}{2} + \sum_{j>i} a_{ij}y_j - b_i.
\end{align}

Our previous work~\cite{ms18} proves the following equivalence result.
\begin{theorem}[\cite{ms18}]
\label{th:cond}
The SOR method and the scheme \eqref{dgscheme1} are equivalent
if
\begin{align}\label{cond1}
h = \frac{2\omega}{2-\omega}. 
\end{align}
In other words, under the relation \eqref{cond1} the sequences generated by
the SOR method and the scheme \eqref{dgscheme1} coincide with each other
as long as the same initial vector $\bx_0$ is used.
\end{theorem}

Therefore, the SOR method can be regarded as an algorithm solving the optimisation problem
for the function \eqref{qf}.

\begin{remark}
It is well known that the SOR method applied to symmetric positive definite linear systems converges
if and only if $\omega \in (0,2)$. 
In the standard SOR theory,
this condition is obtained by analysing the spectral radius of
$\Go$.
However, this condition can also be proved in a completely different manner based on the above discussion.
The convergence condition in terms of $h$ is $h>0$ (as discussed in the previous subsection),
and $h>0$ is equivalent to
$\omega \in (0,2)$ under the relation \eqref{cond1}.
\end{remark}

\section{Adaptive SOR methods}
\label{sec3}

Based on
the interpretation of the SOR method discussed in the previous section,
we develop some new adaptive SOR methods.
The key idea is to control the step size instead of the relaxation parameter.
In the following, the coefficient matrix $A$ is always assumed to be symmetric positive definite.

Note that the step size could be adaptively controlled 
in many different ways.
For example, several step size control techniques have been developed in the context of the numerical
analysis of ordinary differential equations~\cite{hw96}, and thus
applying such techniques would be a possibility~\cite{hw96}.
However, this merely provides an efficient strategy, because
such step size control techniques aim to integrate ODEs as precisely as possible,
but do not aim to find an equilibrium as quickly as possible.

Alternatively, we consider some line search strategies, developed in the context of unconditioned optimisation problems. 
At each iteration, after defining the step direction $\bp^{(k)}$, for example as 
the direction of steepest descent
$\bp_k = -\nabla f(\bx^{(k)})$,
line search methods seek to determine the step size such that it minimises 
$f(\bx + h^{(k)} \bp^{(k)})$.
However, because this task is often computationally expensive for a general function $f$,
most line search methods seek to approximate the optimal step size at a low cost.
Among these, methods employing the Armijo condition or Wolfe conditions are popular.
We hope to apply such line search methods to the SOR method.
However, such line search methods cannot be directly incorporated with the SOR method
because the step direction, i.e., the discrete gradient $\dg f (\bx^{(n+1)},\bx^{(n)})$,
depends on the choice of the step size (note that $\bx^{(n+1)}$ 
depends on the step size).
Therefore, 
slight modifications are required
to apply the aforementioned line search methods to the SOR method.
To achieve this task, we shall directly use the approach proposed in~\cite{rl17}.

In the following, in order to simplify the presentation 
we assume without loss of generality that $D=I$. That is,
the coefficient matrix $A$ is preconditioned 
such that its diagonal matrix coincides with the identity matrix.
This can be achieved by $D^{-1/2} A D^{-1/2}$.

\subsection{Approach based on the locally optimal step size of the steepest descent method}

Our previous work~\cite{ms18} considers an adaptive SOR method
that employs the locally optimal step size of the steepest descent method.
We review this method
before presenting new approaches.

The locally optimal step size of the steepest descent method is obtained by solving
\begin{align}  \label{sorsizeprob}
\min _{h^{(k)} > 0} f(\bx^{(k)} + h^{(k)} \bp^{(k)}), \quad
\bp^{(k)} = - \nabla f(\bx^{(k)}).
\end{align}
For the function \eqref{qf}, the locally optimal step size
is given explicitly by
\begin{align} \label{oph}
h^{(k)} = \frac{(\br^{(k)}) ^\top  \br^{(k)}}{(\br^{(k)})^\top A\br^{(k)}},
\end{align}
where $\br^{(k)} = \bb - A\bx^{(k)}$.
The main idea presented in our previous paper~\cite{ms18}
is to adopt \eqref{oph}
as the step size of the SOR method.
This approach is summarised in Algorithm~\ref{algo:adap1}.

\begin{algorithm}
\caption{Adaptive SOR method based on the steepest descent method}
\begin{algorithmic}
\label{algo:adap1}
\STATE $\br^{(0)} := \bb - A\bx^{(0)}$;
\FOR{$k=0,1,2,\dots$ \textbf{until} $\| \br^{(k)} \| \leq \epsilon \| \bb\|$}
\STATE $\displaystyle h^{(k)} := \frac{(\br^{(k)}) ^\top  \br^{(k)}}{(\br^{(k)})^\top A\br^{(k)}}$;
\STATE $\displaystyle \omega^{(k)} := \frac{2h^{(k)}}{2+h^{(k)}}$;
\STATE $\bx^{(k+1)} := G_{\omega^{(k)}}\bx^{(k)} + \bc$;
%\STATE $\bx^{(k+1)} := (D+\omega^{(k)} L)^{-1} \{[(1-\omega^{(k)})D - \omega^{(k)} U]\bx^{(k)} + \omega^{(k)} \bb \}$;
\STATE $\br^{(k+1)} := \bb - A\bx^{(k+1)}$;
\ENDFOR
\end{algorithmic}
\end{algorithm}

This algorithm does not contain any parameters to be predetermined.
It is observed in~\cite{ms18} that
this algorithm performs well for the Poisson equation:
the number of iterations required for this algorithm is under two times that of
the SOR method with the optimal parameter.
It should be noted that
there is only one additional matrix-vector product per iteration,
i.e., the calculation of
$A \br^{(k)}$ appearing in the denominator of \eqref{oph}. This is
in contrast to some existing approaches, such as~\cite{ba03},
for which additional matrix-vector products are required.

However, it turns out that this approach does not always work perfectly 
when tested with other linear systems,
as illustrated in the next section.
Moreover, it is hoped that the cost for estimating the parameter can be reduced further.

\begin{remark} \label{rem1}
The computational effort for updating the step size (or the relaxation parameter)
could be reduced if the frequency of the update is decreased, and this idea can be applied to
any adaptive SOR method, including those presented in this paper.
\end{remark}

\subsection{New adaptive SOR methods}

Instead of using the step size \eqref{sorsizeprob},
we propose new strategies utilising the Armijo condition or Wolfe conditions
with slight modifications.

Let us start the discussion by applying the Armijo condition,
which is motivated by the recent paper~\cite{rl17}.
As explained in~\cite{rl17},
if $\bx^{(k+1)}$ satisfies
\begin{align} \label{armijo1}
f(\bx^{(k+1)}) \leq f(\bx^{(k)}) + c_1 \nabla f (\bx^{(k)}) ^\top (\bx^{(k+1)} - \bx^{(k)})
\end{align}
for a predetermined constant $c_1 \in (0,1)$,
then the step size is deemed to be good, and increased for the next iteration by a factor of $\lambda_1 >1$:
$h^{(k+1)} = \lambda_1 h^{(k)}$.
If the condition \eqref{armijo1} is not satisfied, then the step size is decreased for the next iteration by a 
factor of $\rho_1 \in(0,1)$: $h^{(k+1)} = \rho_1 h^{(k)}$.
This approach has the following distinguished features.
\begin{itemize}
\item
In contrast to the standard application of the Armijo condition,
the condition \eqref{armijo1} is not used for calculating the current step size $h^{(k)}$,
but rather for calculating the next step size $h^{(k+1)}$.
This is because whether or not the condition is satisfied, the dissipation property
$f(\bx^{(k+1)}) \leq f(\bx^{(k)})$ always follows.
\item 
Assume that the residual $\br^{(k)} = \bb - A\bx^{(k)}$ is calculated and checked at each iteration
(note that this can be obtained without calculating additional matrix-vector products, because
$L\bx^{(k)}$ and $U\bx^{(k)}$ are obtained in the preceding SOR iteration).
Then, there are no additional matrix-vector products required to calculate $f(\bx^{(k)})$
and $\nabla f(\bx^{(k)}) = -\br^{(k)}$.
Thus, the additional costs for checking the condition \eqref{armijo1} consist of only a few
inner product calculations.
\end{itemize}

\begin{remark}
In the condition \eqref{armijo1}, one might think that 
the gradient $\nabla f (\bx^{(k)})$ should be replaced by the 
discrete gradient $\dg f (\bx^{(k+1)},\bx^{(k)})$.
However,
considering the condition
\begin{align} %\label{armijo1}
f(\bx^{(k+1)}) \leq f(\bx^{(k)}) + c_1 \dg f (\bx^{(k+1)},\bx^{(k)}) ^\top (\bx^{(k+1)} - \bx^{(k)})
\end{align}
is meaningless, because this condition is always satisfied:
\begin{align}
f(\bx^{(k+1)}) 
&= f(\bx^{(k)}) +\dg f (\bx^{(k+1)},\bx^{(k)}) ^\top (\bx^{(k+1)} - \bx^{(k)}) \\
&\leq f(\bx^{(k)}) + c_1 \dg f (\bx^{(k+1)},\bx^{(k)}) ^\top (\bx^{(k+1)} - \bx^{(k)}).
\end{align}
The equality follows thanks to \eqref{dcr}, and
the inequality follows from the fact that $f (\bx^{(k+1)},\bx^{(k)}) ^\top (\bx^{(k+1)} - \bx^{(k)})\leq 0$ (see~\eqref{d:dp})
 and the assumption $c_1 \in (0,1)$.
 \end{remark}
 
Our preliminary numerical experiments show that the above approach performs well
for many linear systems.
However, the step size detected using the above approach may become close to zero 
or quite large, which delays the convergence. 
In other words, in terms of the relaxation parameter $\omega$,
$\omega^{(k)}$ sometimes takes a value close to $0$ or $2$.
To avoid this situation, we need to set the maximum and minimum values 
($M_\omega$ and $\epsilon_\omega$) that the relaxation parameter is allowed to take,
and if this is violated 
we reset the relaxation parameter and the step size.

The above algorithm is summarised in Algorithm~\ref{algo:adap2}.
Here, the step size is initially set to $2$, which indicates that the relaxation parameter
is initially set to $1$, and if $\omega^{(k+1)} \notin (\epsilon_\omega , M_\omega)$ then 
the step size and relaxation parameter are reset to $1$ and $2$, respectively.

\begin{algorithm}
\caption{Adaptive SOR method based on the Armijo condition}
\begin{algorithmic}
\label{algo:adap2}
\STATE Set parameters $c_1 \in (0,1)$, $\lambda_1 > 1$,
$\rho_1 \in (0,1)$;
\STATE $\br^{(0)} := \bb - A\bx^{(0)}$;
\STATE $h^{(0)} := 2$;
\STATE $\omega^{(0)}:=1$;
\FOR{$k=0,1,2,\dots$ \textbf{until} $\| \br^{(k)} \| \leq \epsilon \| \bb\|$}
\STATE $\bx^{(k+1)} := G_{\omega^{(k)}}\bx^{(k)} + \bc$;
%\STATE $\bx^{(k+1)} := (D+\omega^{(k)} L)^{-1} \{[(1-\omega^{(k)})D - \omega^{(k)} U]\bx^{(k)} + \omega^{(k)} \bb \}$;
\STATE $\br^{(k+1)} := \bb - A\bx^{(k+1)}$;
\IF{$f(\bx^{(k+1)}) \leq f(\bx^{(k)}) + c_1 \nabla f (\bx^{(k)}) ^\top (\bx^{(k+1)} - \bx^{(k)})$}
\STATE $h^{(k+1)} := \lambda_1 h^{(k)}$;
\ELSE
\STATE $h^{(k+1)} := \rho_1 h^{(k)}$
\ENDIF
\STATE $\displaystyle \omega^{(k+1)} := \frac{2h^{(k+1)}}{2+h^{(k+1)}}$;
\IF{$\omega^{(k+1)} \notin (\epsilon_\omega , M_\omega)$}
\STATE  $h^{(k+1)} = 2$;
\STATE  $\omega^{(k+1)} = 1$;
\ENDIF
\ENDFOR
\end{algorithmic}
\end{algorithm}

In the context of nonlinear optimisation problems, 
the Wolfe conditions, which consist of the Armijo condition and the so-called curvature condition,
is often employed to avoid the step size becoming close to $0$.
In the following, we consider another algorithm, which is based on
the Wolfe conditions.

Here, we consider the following curvature condition:
\begin{align}\label{wolfe1}
c_2 \nabla f (\bx^{(k)}) ^\top (\bx^{(k+1)} - \bx^{(k)} )\leq \nabla f(\bx^{(k+1)})^\top (\bx^{(k+1)}-\bx^{(k)})
\end{align}
for a predetermined constant $c_2\in (c_1,1)$.
If the two conditions \eqref{armijo1} and \eqref{wolfe1} are satisfied, 
then the step size for the next iteration is increased by $h^{(k+1)} = \lambda_1 h^{(k)}$, as in Algorithm~1.
However, if the condition \eqref{wolfe1} is not satisfied, then the current step size is deemed too small
and increased for the next iteration by a larger factor of $\lambda_2>\lambda_1 (>1)$:
$h^{(k+1)} = \lambda_2 h^{(k)}$.

The above approach is summarised in Algorithm~\ref{algo:adap3}.

\begin{algorithm}
\caption{Adaptive SOR method based on the Wolfe conditions}
\begin{algorithmic}
\label{algo:adap3}
\STATE Set parameters $c_1 \in (0,1)$, $c_2 \in (c_1, 1)$, $\lambda_1 > 1$,
$\lambda_2 > \lambda_1$, $\rho_1 \in (0,1)$;
\STATE $\br^{(0)} := \bb - A\bx^{(0)}$;
\STATE $h^{(0)} := 2$;
\STATE $\omega^{(0)}:=1$;
\FOR{$k=0,1,2,\dots$ \textbf{until} $\| \br^{(k)} \| \leq \epsilon \| \bb\|$}
\STATE $\bx^{(k+1)} := G_{\omega^{(k)}}\bx^{(k)} + \bc$;
\STATE $\br^{(k+1)} := \bb - A\bx^{(k+1)}$;
\IF{$f(\bx^{(k+1)}) \leq f(\bx^{(k)}) + c_1 \nabla f (\bx^{(k)}) ^\top (\bx^{(k+1)} - \bx^{(k)})$}
\IF{$c_2 \nabla f (\bx^{(k)}) ^\top (\bx^{(k+1)} - \bx^{(k)} )\leq \nabla f(\bx^{(k+1)})^\top (\bx^{(k+1)}-\bx^{(k)})$}
\STATE $h^{(k+1)} := \lambda_1 h^{(k)}$;
\ELSE
\STATE $h^{(k+1)} := \lambda_2 h^{(k)}$;
\ENDIF
\ELSE
\STATE $h^{(k+1)} := \rho_1 h^{(k)}$
\ENDIF
\STATE $\displaystyle \omega^{(k+1)} := \frac{2h^{(k+1)}}{2+h^{(k+1)}}$;
\IF{$\omega^{(k+1)} \notin (\epsilon_\omega , M_\omega)$}
\STATE  $h^{(k+1)} = 2$;
\STATE  $\omega^{(k+1)} = 1$;
\ENDIF
\ENDFOR
\end{algorithmic}
\end{algorithm}

\section{Numerical experiments}
\label{sec4}

In this section, the proposed adaptive SOR methods are tested
for a variety of symmetric positive definite linear systems.
The main aim of this section is 
to numerically verify the convergence of the proposed methods.
Note that the computational costs required for updating the step size 
were discussed from the viewpoint of the number of matrix-vector products 
in the previous section, 
and the actual computation time depends on 
the frequency of the updates, as explained in Remark~\ref{rem1}.
Therefore, we do not check the actual computational time, and in all
numerical experiments the step size is updated at every iteration, for a
fair comparison with Algorithms~\ref{algo:adap1},~\ref{algo:adap2}, 
and \ref{algo:adap3}.
We also note that a comparison with other types of adaptive SOR methods or linear solvers, 
such as the conjugate gradient method, is beyond the scope of this study,
because of the difficulty in setting a
fair gauge for a head-to-head comparison owing to quite different mathematical
features.
All numerical experiments are conducted 
after the coefficient matrix $A$ is preconditioned 
such that its diagonal matrix coincides with the identity matrix.
All the computations are performed in a computation environment with 
3.5 GHz Intel Core i5, 8 GB memory, and OS X 10.10.5. We employ MATLAB (R2015a).
Below, the
results are often displayed after thinning out the data to highlight the behaviour of each algorithm and reduce
the file size.

Note that in Algorithms~\ref{algo:adap2} and~\ref{algo:adap3} 
there are several constants and factors to be predetermined.
Instead of finding optimal values for each linear system, 
which seems more difficult than
determining the optimal relaxation parameter of the SOR method, 
we fix them as $(c_1,c_2,\lambda_1,\lambda_2,\rho_1)
= (0.89, 0.95,1.15,1.4,0.85 )$,
and this combination is employed for different linear systems.

\subsection{Poisson equation}
The first problem arises from the finite difference discretisation of 
the Poisson equation
\begin{align}
-\triangle u  &= f \quad \text{in }\Omega, \\
u&=0 \quad \text{on }\partial \Omega,
\end{align}
where $\Omega = (0,1)^2\subset \bbR^2$, $\triangle$ 
denotes the Laplace operator
and $f:\Omega \to \bbR$.
The standard finite difference discretisation on a uniform mesh $N\times N$
leads to a symmetric positive definite linear system
for which the coefficient matrix $A$ is given by
\begin{align}
A = \frac{1}{\Delta x^2} \left(-I \otimes I + \frac{1}{4} B\otimes I + \frac{1}{4} I \otimes B \right) ,
\end{align}
where $I\in\bbR^{(N-1)\times (N-1)}$ is the identity matrix,
\begin{align*}
B = \begin{bmatrix}
0 & 1 & &  \\
1 & \ddots & \ddots& \\
& \ddots & \ddots & 1 \\
& &1 & 0 
\end{bmatrix} \in \bbR^{(N-1)\times (N-1)},
\end{align*}
and $\otimes$ denotes the Kronecker product.
The optimal relaxation parameter for this matrix $A$ is well known, 
and is expressed in terms of $h=1/(N+1)$ as
\begin{align*}
\omega_{\text{opt}} = \frac{2}{1+\sqrt{1-\cos^2 (\pi h)}}.
\end{align*}

Fig.~\ref{fig:poisson} illustrates the numerical results 
for a test problem $f(x,y) = \sin(\pi x) \sin (\pi y)$.
The left figures show the relative residual 2-norm. 
The convergence of Algorithm~\ref{algo:adap1}
is slower than that of the SOR method with the optimal parameter, 
but is considerably faster than that of the Gauss-Seidel method. Furthermore,
the number of iterations required for convergence is less than two times 
that of the SOR method with the optimal parameter.
In this sense, as discussed in our previous report~\cite{ms18}, 
Algorithm~\ref{algo:adap1} appears to be flexible.
However, as previously explained, 
the computational cost per iteration of Algorithm~\ref{algo:adap1}
is almost twice as expensive as that of the SOR method.
Keeping this in mind, 
let us discuss the behaviour of Algorithms~\ref{algo:adap2} and~\ref{algo:adap3}.
For $N=60$, both algorithms perform satisfactorily compared 
with Algorithm~\ref{algo:adap1}.
However, for $N=100$ and $N=120$ there are substantial differences.
While the results for Algorithm~\ref{algo:adap3} still remain satisfactory, the 
convergence of Algorithm~\ref{algo:adap2} appears to deteriorate.
These differences can also be explained from the right-side figures, 
which show the variations of $\omega^{(k)}$.
In particular, the difference is remarkable when $N=100$. 
Here, $\omega^{(k)}$ of Algorithm~\ref{algo:adap2}
tends to $0$, even after it is reset to $1$ at around the $1500$th iteration.

\input{fig_poisson.tex}

\subsection{Additional examples}

Table~\ref{tab:matrix_data} presents the matrix data used in 
the numerical experiments,
which were obtained from 
the Matrix Market (\url{https://math.nist.gov/MatrixMarket/}) 
and the ELSES Matrix Library (\url{http://www.elses.jp/matrix/}).
In the following numerical experiments, 
the right-hand vector was set to $\bb = (1,\dots,1)^\top$.

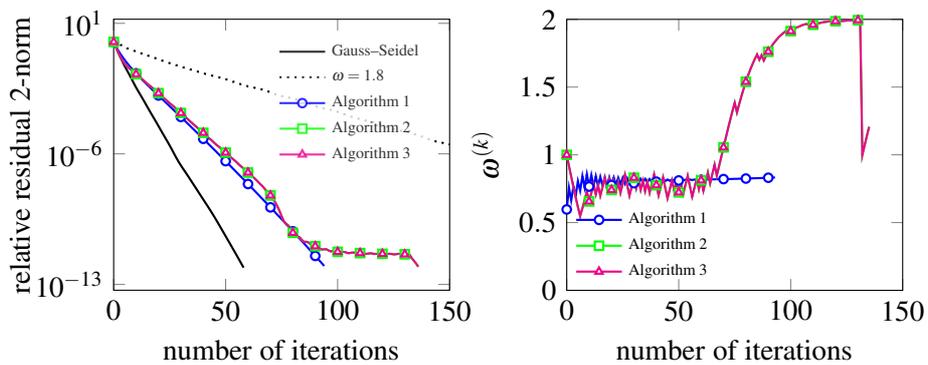
\begin{figure}[htbp]
\centering

\begin{tabular}{ll}
\begin{tikzpicture}
\tikzstyle{every node}=[]
\begin{semilogyaxis}[width=6cm,%restrict x to domain = 0:69,
xmax=150,xmin=0,
%	restrict y to domain=-0.8:1.1, ymax=1.1,ymin=-0.8,
xlabel={number of iterations},ylabel={relative residual $2$-norm},
legend entries={Gauss--Seidel,$\omega=1.8$,Algorithm~\ref{algo:adap1},Algorithm~\ref{algo:adap2},Algorithm~\ref{algo:adap3}},
legend style={font = \tiny ,legend cell align=left,draw=none,fill=white, fill opacity=0.8, text opacity=1,at={(.95,.95)},anchor=north east}
	]
%\node[above right] at (axis cs: -2,-1) {$h=0.1$};
%\node[above left] at (axis cs: 0.5,-1.05) {$60$ steps};
\addplot[thick] table {
0 1.00e+00
1 5.39e-01
2 2.83e-01
3 1.52e-01
4 8.41e-02
5 4.87e-02
6 2.93e-02
7 1.76e-02
8 1.04e-02
9 6.08e-03
10 3.62e-03
11 2.24e-03
12 1.41e-03
13 8.82e-04
14 5.44e-04
15 3.33e-04
16 2.05e-04
17 1.28e-04
18 8.03e-05
19 5.02e-05
20 3.10e-05
21 1.91e-05
22 1.18e-05
23 7.35e-06
24 4.58e-06
25 2.82e-06
26 1.70e-06
27 1.01e-06
28 6.07e-07
29 3.75e-07
30 2.41e-07
31 1.59e-07
32 1.06e-07
33 7.02e-08
34 4.61e-08
35 3.01e-08
36 1.97e-08
37 1.29e-08
38 8.41e-09
39 5.50e-09
40 3.61e-09
41 2.36e-09
42 1.55e-09
43 1.00e-09
44 6.47e-10
45 4.12e-10
46 2.60e-10
47 1.62e-10
48 1.01e-10
49 6.25e-11
50 3.87e-11
51 2.40e-11
52 1.49e-11
53 9.21e-12
54 5.70e-12
55 3.51e-12
56 2.15e-12
57 1.31e-12
58 7.98e-13
};
\addplot[thick,dotted] table {
0 1.00e+00
1 1.02e+00
2 9.29e-01
3 6.99e-01
4 6.57e-01
5 5.75e-01
6 4.56e-01
7 4.15e-01
8 3.94e-01
9 3.40e-01
10 3.13e-01
11 2.72e-01
12 2.64e-01
13 2.16e-01
14 2.02e-01
15 1.98e-01
16 1.68e-01
17 1.52e-01
18 1.47e-01
19 1.29e-01
20 1.13e-01
21 1.06e-01
22 9.69e-02
23 8.89e-02
24 7.90e-02
25 7.34e-02
26 6.93e-02
27 6.27e-02
28 5.63e-02
29 5.23e-02
30 5.00e-02
31 4.38e-02
32 4.02e-02
33 3.91e-02
34 3.44e-02
35 3.18e-02
36 3.09e-02
37 2.71e-02
38 2.48e-02
39 2.42e-02
40 2.17e-02
41 1.92e-02
42 1.83e-02
43 1.73e-02
44 1.56e-02
45 1.39e-02
46 1.32e-02
47 1.28e-02
48 1.13e-02
49 9.98e-03
50 9.95e-03
51 9.46e-03
52 8.02e-03
53 7.34e-03
54 7.43e-03
55 6.76e-03
56 5.69e-03
57 5.42e-03
58 5.40e-03
59 4.80e-03
60 4.09e-03
61 3.97e-03
62 3.93e-03
63 3.45e-03
64 2.96e-03
65 2.93e-03
66 2.89e-03
67 2.49e-03
68 2.15e-03
69 2.19e-03
70 2.13e-03
71 1.78e-03
72 1.58e-03
73 1.64e-03
74 1.55e-03
75 1.27e-03
76 1.18e-03
77 1.23e-03
78 1.11e-03
79 9.06e-04
80 8.82e-04
81 9.01e-04
82 7.81e-04
83 6.51e-04
84 6.53e-04
85 6.47e-04
86 5.46e-04
87 4.67e-04
88 4.79e-04
89 4.60e-04
90 3.79e-04
91 3.37e-04
92 3.50e-04
93 3.24e-04
94 2.63e-04
95 2.47e-04
96 2.54e-04
97 2.25e-04
98 1.86e-04
99 1.83e-04
100 1.81e-04
101 1.56e-04
102 1.35e-04
103 1.34e-04
104 1.28e-04
105 1.11e-04
106 9.86e-05
107 9.68e-05
108 9.18e-05
109 7.99e-05
110 7.16e-05
111 7.04e-05
112 6.65e-05
113 5.75e-05
114 5.21e-05
115 5.16e-05
116 4.81e-05
117 4.12e-05
118 3.80e-05
119 3.79e-05
120 3.47e-05
121 2.96e-05
122 2.79e-05
123 2.78e-05
124 2.49e-05
125 2.12e-05
126 2.06e-05
127 2.04e-05
128 1.78e-05
129 1.52e-05
130 1.53e-05
131 1.49e-05
132 1.26e-05
133 1.11e-05
134 1.14e-05
135 1.08e-05
136 8.88e-06
137 8.23e-06
138 8.49e-06
139 7.65e-06
140 6.33e-06
141 6.15e-06
142 6.22e-06
143 5.41e-06
144 4.58e-06
145 4.60e-06
146 4.50e-06
147 3.82e-06
148 3.36e-06
149 3.41e-06
150 3.22e-06
};
\addplot[thick,mark=*,mark repeat={10},mark options={scale=0.8,fill=white},color=blue] table {
0 1.00e+00
1 5.76e-01
2 3.54e-01
3 2.28e-01
4 1.54e-01
5 1.05e-01
6 7.41e-02
7 5.21e-02
8 3.84e-02
9 2.76e-02
10 2.09e-02
11 1.53e-02
12 1.17e-02
13 8.68e-03
14 6.71e-03
15 5.02e-03
16 3.91e-03
17 2.94e-03
18 2.29e-03
19 1.73e-03
20 1.35e-03
21 1.02e-03
22 7.97e-04
23 6.04e-04
24 4.70e-04
25 3.57e-04
26 2.77e-04
27 2.10e-04
28 1.63e-04
29 1.24e-04
30 9.56e-05
31 7.25e-05
32 5.60e-05
33 4.25e-05
34 3.27e-05
35 2.48e-05
36 1.91e-05
37 1.45e-05
38 1.11e-05
39 8.41e-06
40 6.43e-06
41 4.87e-06
42 3.72e-06
43 2.82e-06
44 2.15e-06
45 1.63e-06
46 1.24e-06
47 9.36e-07
48 7.11e-07
49 5.38e-07
50 4.08e-07
51 3.08e-07
52 2.33e-07
53 1.76e-07
54 1.33e-07
55 1.00e-07
56 7.57e-08
57 5.70e-08
58 4.30e-08
59 3.23e-08
60 2.43e-08
61 1.83e-08
62 1.37e-08
63 1.03e-08
64 7.74e-09
65 5.80e-09
66 4.35e-09
67 3.26e-09
68 2.44e-09
69 1.82e-09
70 1.36e-09
71 1.02e-09
72 7.59e-10
73 5.66e-10
74 4.22e-10
75 3.14e-10
76 2.33e-10
77 1.74e-10
78 1.29e-10
79 9.57e-11
80 7.10e-11
81 5.26e-11
82 3.89e-11
83 2.88e-11
84 2.13e-11
85 1.57e-11
86 1.16e-11
87 8.53e-12
88 6.28e-12
89 4.62e-12
90 3.39e-12
91 2.49e-12
92 1.83e-12
93 1.34e-12
94 9.78e-13
};
\addplot[thick,mark=square*,mark repeat={10},mark options={scale=0.8,fill=white},color=green] table {
0 1.00e+00
1 5.39e-01
2 2.75e-01
3 1.54e-01
4 9.51e-02
5 6.50e-02
6 4.82e-02
7 3.80e-02
8 3.02e-02
9 2.39e-02
10 1.88e-02
11 1.48e-02
12 1.18e-02
13 9.40e-03
14 7.32e-03
15 5.86e-03
16 4.60e-03
17 3.71e-03
18 2.95e-03
19 2.26e-03
20 1.79e-03
21 1.37e-03
22 1.09e-03
23 8.43e-04
24 6.72e-04
25 5.23e-04
26 4.19e-04
27 3.30e-04
28 2.49e-04
29 1.94e-04
30 1.56e-04
31 1.24e-04
32 9.47e-05
33 7.48e-05
34 5.73e-05
35 4.54e-05
36 3.50e-05
37 2.78e-05
38 2.16e-05
39 1.72e-05
40 1.35e-05
41 1.01e-05
42 7.89e-06
43 6.31e-06
44 4.98e-06
45 3.79e-06
46 2.98e-06
47 2.40e-06
48 1.91e-06
49 1.47e-06
50 1.17e-06
51 9.05e-07
52 7.20e-07
53 5.60e-07
54 4.47e-07
55 3.51e-07
56 2.65e-07
57 2.07e-07
58 1.65e-07
59 1.31e-07
60 9.97e-08
61 7.83e-08
62 6.01e-08
63 4.73e-08
64 3.66e-08
65 2.71e-08
66 2.08e-08
67 1.56e-08
68 1.10e-08
69 8.16e-09
70 5.89e-09
71 4.08e-09
72 2.68e-09
73 1.64e-09
74 9.04e-10
75 4.48e-10
76 2.41e-10
77 1.79e-10
78 1.43e-10
79 9.94e-11
80 5.91e-11
81 4.06e-11
82 3.87e-11
83 3.38e-11
84 2.45e-11
85 1.90e-11
86 1.95e-11
87 1.78e-11
88 1.40e-11
89 1.12e-11
90 1.12e-11
91 1.06e-11
92 8.71e-12
93 7.41e-12
94 7.71e-12
95 7.69e-12
96 6.59e-12
97 5.90e-12
98 6.32e-12
99 6.38e-12
100 5.60e-12
101 5.19e-12
102 5.50e-12
103 5.55e-12
104 5.06e-12
105 4.75e-12
106 5.03e-12
107 5.17e-12
108 4.78e-12
109 4.51e-12
110 4.81e-12
111 4.90e-12
112 4.53e-12
113 4.37e-12
114 4.61e-12
115 4.65e-12
116 4.42e-12
117 4.32e-12
118 4.43e-12
119 4.51e-12
120 4.39e-12
121 4.22e-12
122 4.25e-12
123 4.42e-12
124 4.33e-12
125 4.08e-12
126 4.16e-12
127 4.40e-12
128 4.27e-12
129 4.00e-12
130 4.13e-12
131 4.36e-12
132 4.20e-12
133 2.68e-12
134 1.90e-12
135 1.34e-12
136 8.93e-13
};
\addplot[thick,mark=triangle*,mark repeat={10},mark options={scale=0.8,fill=white},color=magenta] table {
0 1.00e+00
1 5.39e-01
2 2.75e-01
3 1.54e-01
4 9.51e-02
5 6.50e-02
6 4.82e-02
7 3.80e-02
8 3.02e-02
9 2.39e-02
10 1.88e-02
11 1.48e-02
12 1.18e-02
13 9.40e-03
14 7.32e-03
15 5.86e-03
16 4.60e-03
17 3.71e-03
18 2.95e-03
19 2.26e-03
20 1.79e-03
21 1.37e-03
22 1.09e-03
23 8.43e-04
24 6.72e-04
25 5.23e-04
26 4.19e-04
27 3.30e-04
28 2.49e-04
29 1.94e-04
30 1.56e-04
31 1.24e-04
32 9.47e-05
33 7.48e-05
34 5.73e-05
35 4.54e-05
36 3.50e-05
37 2.78e-05
38 2.16e-05
39 1.72e-05
40 1.35e-05
41 1.01e-05
42 7.89e-06
43 6.31e-06
44 4.98e-06
45 3.79e-06
46 2.98e-06
47 2.40e-06
48 1.91e-06
49 1.47e-06
50 1.17e-06
51 9.05e-07
52 7.20e-07
53 5.60e-07
54 4.47e-07
55 3.51e-07
56 2.65e-07
57 2.07e-07
58 1.65e-07
59 1.31e-07
60 9.97e-08
61 7.83e-08
62 6.01e-08
63 4.73e-08
64 3.66e-08
65 2.71e-08
66 2.08e-08
67 1.56e-08
68 1.10e-08
69 8.16e-09
70 5.89e-09
71 4.08e-09
72 2.68e-09
73 1.64e-09
74 9.04e-10
75 4.48e-10
76 2.41e-10
77 1.79e-10
78 1.43e-10
79 9.94e-11
80 5.91e-11
81 4.06e-11
82 3.87e-11
83 3.38e-11
84 2.45e-11
85 1.90e-11
86 1.95e-11
87 1.78e-11
88 1.40e-11
89 1.12e-11
90 1.12e-11
91 1.06e-11
92 8.71e-12
93 7.41e-12
94 7.71e-12
95 7.69e-12
96 6.59e-12
97 5.90e-12
98 6.32e-12
99 6.38e-12
100 5.60e-12
101 5.19e-12
102 5.50e-12
103 5.55e-12
104 5.06e-12
105 4.75e-12
106 5.03e-12
107 5.17e-12
108 4.78e-12
109 4.51e-12
110 4.81e-12
111 4.90e-12
112 4.53e-12
113 4.37e-12
114 4.61e-12
115 4.65e-12
116 4.42e-12
117 4.32e-12
118 4.43e-12
119 4.51e-12
120 4.39e-12
121 4.22e-12
122 4.25e-12
123 4.42e-12
124 4.33e-12
125 4.08e-12
126 4.16e-12
127 4.40e-12
128 4.27e-12
129 4.00e-12
130 4.13e-12
131 4.36e-12
132 4.20e-12
133 2.68e-12
134 1.90e-12
135 1.34e-12
136 8.93e-13
};
\end{semilogyaxis}
\end{tikzpicture}
&\hspace{-2em}
\begin{tikzpicture}
\tikzstyle{every node}=[]
\begin{axis}[width=6cm,%restrict x to domain = 0:69,
xmax=150,xmin=0,
%restrict y to domain=0:2,
ymax=2,ymin=0,
xlabel={number of iterations},ylabel={$\omega^{(k)}$},
ylabel near ticks,
legend entries={Algorithm~\ref{algo:adap1},Algorithm~\ref{algo:adap2},Algorithm~\ref{algo:adap3}},
legend style={font = \scriptsize,legend pos=south west,legend cell align=left,draw=none,fill=none},
legend style={font = \tiny ,legend cell align=left,draw=none,fill=none,at={(0,0)},anchor=south west}
	]
\addplot[thick,mark=*,mark repeat={10},mark options={scale=0.8,fill=white},color=blue] table {
0 0.597
1 0.758
2 0.663
3 0.767
4 0.698
5 0.806
6 0.735
7 0.844
8 0.759
9 0.857
10 0.766
11 0.851
12 0.767
13 0.843
14 0.769
15 0.838
16 0.771
17 0.836
18 0.774
19 0.834
20 0.777
21 0.833
22 0.780
23 0.831
24 0.783
25 0.830
26 0.786
27 0.828
28 0.788
29 0.827
30 0.791
31 0.826
32 0.793
33 0.824
34 0.796
35 0.823
36 0.798
37 0.822
38 0.800
39 0.821
40 0.802
41 0.820
42 0.804
43 0.819
44 0.806
45 0.818
46 0.808
47 0.818
48 0.809
49 0.818
50 0.811
51 0.818
52 0.812
53 0.817
54 0.813
55 0.818
56 0.814
57 0.818
58 0.815
59 0.818
60 0.816
61 0.818
62 0.817
63 0.819
64 0.818
65 0.819
66 0.819
67 0.820
68 0.820
69 0.821
70 0.821
71 0.821
72 0.821
73 0.822
74 0.822
75 0.823
76 0.823
77 0.824
78 0.824
79 0.825
80 0.825
81 0.826
82 0.826
83 0.827
84 0.827
85 0.828
86 0.828
87 0.829
88 0.830
89 0.830
90 0.831
91 0.832
92 0.832
93 0.833
};
\addplot[thick,mark=square*,mark repeat={10},mark options={scale=0.8,fill=white},color=green] table {
0 1.000
1 0.919
2 0.839
3 0.761
4 0.686
5 0.615
6 0.548
7 0.605
8 0.666
9 0.729
10 0.656
11 0.719
12 0.784
13 0.708
14 0.773
15 0.698
16 0.762
17 0.829
18 0.752
19 0.818
20 0.741
21 0.807
22 0.731
23 0.796
24 0.720
25 0.786
26 0.853
27 0.775
28 0.699
29 0.764
30 0.831
31 0.753
32 0.820
33 0.743
34 0.809
35 0.732
36 0.798
37 0.721
38 0.787
39 0.855
40 0.776
41 0.701
42 0.765
43 0.833
44 0.755
45 0.680
46 0.744
47 0.810
48 0.734
49 0.800
50 0.723
51 0.789
52 0.712
53 0.778
54 0.845
55 0.767
56 0.692
57 0.756
58 0.823
59 0.746
60 0.812
61 0.735
62 0.801
63 0.869
64 0.790
65 0.858
66 0.927
67 0.847
68 0.916
69 0.985
70 1.055
71 1.124
72 1.193
73 1.259
74 1.323
75 1.384
76 1.313
77 1.374
78 1.433
79 1.488
80 1.539
81 1.587
82 1.631
83 1.671
84 1.708
85 1.741
86 1.702
87 1.736
88 1.696
89 1.730
90 1.761
91 1.789
92 1.814
93 1.836
94 1.856
95 1.874
96 1.889
97 1.903
98 1.915
99 1.901
100 1.913
101 1.924
102 1.934
103 1.942
104 1.949
105 1.956
106 1.962
107 1.966
108 1.961
109 1.966
110 1.960
111 1.965
112 1.969
113 1.973
114 1.977
115 1.980
116 1.982
117 1.985
118 1.987
119 1.984
120 1.986
121 1.984
122 1.986
123 1.988
124 1.989
125 1.991
126 1.992
127 1.993
128 1.992
129 1.993
130 1.994
131 1.995
132 1.000
133 1.070
134 1.139
135 1.207
};
\addplot[thick,mark=triangle*,mark repeat={10},mark options={scale=0.8,fill=white},color=magenta] table {
0 1.000
1 0.919
2 0.839
3 0.761
4 0.686
5 0.615
6 0.548
7 0.605
8 0.666
9 0.729
10 0.656
11 0.719
12 0.784
13 0.708
14 0.773
15 0.698
16 0.762
17 0.829
18 0.752
19 0.818
20 0.741
21 0.807
22 0.731
23 0.796
24 0.720
25 0.786
26 0.853
27 0.775
28 0.699
29 0.764
30 0.831
31 0.753
32 0.820
33 0.743
34 0.809
35 0.732
36 0.798
37 0.721
38 0.787
39 0.855
40 0.776
41 0.701
42 0.765
43 0.833
44 0.755
45 0.680
46 0.744
47 0.810
48 0.734
49 0.800
50 0.723
51 0.789
52 0.712
53 0.778
54 0.845
55 0.767
56 0.692
57 0.756
58 0.823
59 0.746
60 0.812
61 0.735
62 0.801
63 0.869
64 0.790
65 0.858
66 0.927
67 0.847
68 0.916
69 0.985
70 1.055
71 1.124
72 1.193
73 1.259
74 1.323
75 1.384
76 1.313
77 1.374
78 1.433
79 1.488
80 1.539
81 1.587
82 1.631
83 1.671
84 1.708
85 1.741
86 1.702
87 1.736
88 1.696
89 1.730
90 1.761
91 1.789
92 1.814
93 1.836
94 1.856
95 1.874
96 1.889
97 1.903
98 1.915
99 1.901
100 1.913
101 1.924
102 1.934
103 1.942
104 1.949
105 1.956
106 1.962
107 1.966
108 1.961
109 1.966
110 1.960
111 1.965
112 1.969
113 1.973
114 1.977
115 1.980
116 1.982
117 1.985
118 1.987
119 1.984
120 1.986
121 1.984
122 1.986
123 1.988
124 1.989
125 1.991
126 1.992
127 1.993
128 1.992
129 1.993
130 1.994
131 1.995
132 1.000
133 1.070
134 1.139
135 1.207
};
\end{axis}

\end{tikzpicture}
\end{tabular}

\caption{Numerical results for the test problem for ICNT1800: (LEFT) relative residual 2-norm, and (RIGHT) variations of $\omega^{(k)}$.
In the left figure, the results by the Gauss--Seidel method 
is shown in not the dashed line but the black one, and the results for $\omega = 1.9$ are not displayed
since the Gauss--Seidel method was the best.
}
\label{fig:icnt1800}

\end{figure}

The results are illustrated in Figs.~\ref{fig:bcsstk04},~\ref{fig:bcsstk05},~\ref{fig:bcsstk07}, and~\ref{fig:icnt1800}.
In the left-hand figures, in addition to Algorithms~\ref{algo:adap1},~\ref{algo:adap2}, and~\ref{algo:adap3} 
we plot the results for the Gauss--Seidel method (dashed line), 
the SOR method with $\omega=1.8$ (dotted line),
and the SOR method, for which the convergence is the fastest in $\omega = 0.1,0.2,\dots, 1.9$
(black line).
In Fig.~\ref{fig:icnt1800}, the results for the Gauss--Seidel method 
are illustrated by the solid black line rather than the dashed one,
because the Gauss--Seidel method performed best. 

Although the results differ for each problem, 
Algorithm~\ref{algo:adap3} is faster than 
Algorithms~\ref{algo:adap1} and~\ref{algo:adap2} under all settings.
Note that from Fig.~\ref{fig:icnt1800} one might infer that 
Algorithm~\ref{algo:adap1} is the fastest,
but if we take the computational cost at each iteration into account 
Algorithms~\ref{algo:adap2} and~\ref{algo:adap3} are superior.
The number of iterations required to converge for Algorithm~\ref{algo:adap3} is less than
three times as big as that for the SOR method, for which the convergence 
is the fastest in $\omega = 0.1,0.2,\dots, 1.9$. 
Furthermore, in the results shown in Fig.~\ref{fig:bcsstk07}, Algorithm~\ref{algo:adap3} 
is the fastest.
We observed similar behaviours for the other problems 
obtained from the Matrix Market and the ELSES Matrix Library.

\begin{table}
\centering
\caption{Matrix data.}
\label{tab:matrix_data}
\begin{tabular}{c|c|c}
\hline 
 Data & $n$ & {\#}nz\\ 
\hline 
 BCSSTK04 & 132 & 3648 \\ 
 \hline 
 BCSSTK05 & 153 & 2423 \\ 
 \hline 
 BCSSTK07 & 420 & 7860 \\ 
 \hline 
 ICNT1800 & 1800 & 574960 \\ 
\end{tabular} 
\end{table}

\input{fig_bcsstk04.tex}
\input{fig_bcsstk05.tex}
\input{fig_bcsstk07.tex}

\section{Concluding remarks}
\label{sec5}

In this paper, we proposed two adaptive SOR methods based on line search techniques.
One is based on the Armijo condition, and the other on the Wolfe conditions.
The features of the proposed methods are summarised as follows:
\begin{itemize}
\item They are applicable to any symmetric positive definite linear system,
in the sense that the convergence condition is always satisfied.
\item No additional matrix-vector products are required to update the relaxation parameter.
\item There are several factors and parameters in the algorithms.
While these should be predetermined and the optimal combination is difficult to obtain,
for the adaptive SOR method based on the Wolfe conditions (Algorithm~\ref{algo:adap3}),
the empirical combination  $(c_1,c_2,\lambda_1,\lambda_2,\rho_1)
= (0.89, 0.95,1.15,1.4,0.85 )$ 
employed in this paper can perform well for many symmetric positive definite
linear systems.
\item The convergence is slower than that of the SOR method using the optimal parameter,
but is faster than naive choices such as $\omega = 1$ (the Gauss--Seldel method) and $\omega=1.8$ in most cases.
Furthermore, the number of iterations required for convergence is less than two times that of  
the SOR method with the optimal parameter in most cases,
which indicates that the computational cost is also less than double that of 
the SOR method with the optimal parameter.
\end{itemize}

We mention several directions for future work.
In this paper, among the infinite number of combinations of parameters that are to be 
predetermined, 
the results obtained by the chosen empirical combination were presented, but 
it would be interesting to study a strategy of tuning the combination 
to find broad applications.
We are also currently attempting to extend the presented approach to more general linear systems.

\bibliographystyle{plain}
\bibliography{references}
\end{document}